\documentclass[11pt,a4paper]{amsart}

\usepackage[T1]{fontenc}
\usepackage{lmodern}
\usepackage[UKenglish]{babel}
\usepackage[top=2cm,bottom=2cm,left=2cm,right=2cm]{geometry}
\usepackage{amsmath,amssymb,amsthm,mathtools}
\usepackage{enumitem}
\usepackage{microtype}
\usepackage[numbers,square,sort&compress]{natbib}

\usepackage[colorlinks,pagebackref,hypertexnames=false]{hyperref}

\newcommand{\SO}{\operatorname{SO}}
\newcommand{\GL}{\operatorname{GL}}
\newcommand{\U}{\operatorname{U}}
\newcommand{\Spin}{\operatorname{Spin}}
\newcommand{\Gtwo}{\mathrm{G}_2}
\newcommand{\fh}{\mathfrak{h}}
\newcommand{\fm}{\mathfrak{m}}
\newcommand{\fso}{\mathfrak{so}}
\newcommand{\SU}{\operatorname{SU}}

\newcommand{\Fr}{\operatorname{Fr}}
\newcommand{\End}{\operatorname{End}}
\newcommand{\Stab}{\operatorname{Stab}}
\newcommand{\vol}{\operatorname{vol}}
\newcommand{\Ric}{\operatorname{Ric}}
\newcommand{\id}{\operatorname{id}}
\newcommand{\cC}{\mathcal{C}}
\newcommand{\cQ}{\mathcal{Q}}

\newtheorem{theorem}{Theorem}[section]

\newtheorem{corollary}[theorem]{Corollary}

\theoremstyle{definition}

\theoremstyle{remark}
\newtheorem{remark}[theorem]{Remark}

\numberwithin{equation}{section}

\title[Taylor formula for H-structures]
{Taylor formula for \texorpdfstring{$H$}{H}-structures}

\author[Lucca Delboni]{Lucca Delboni\protect\textsuperscript{\normalfont *}}

\author[Henrique N. Sá Earp]{Henrique N. Sá Earp}

\makeatletter
\g@addto@macro\@setauthors{\begingroup\footnotesize
\begin{center}
Universidade Estadual de Campinas (Unicamp).

\smallskip
Email: \href{mailto:lucca.delboni@ime.unicamp.br}{\nolinkurl{lucca.delboni@ime.unicamp.br}} ;\\
\href{mailto:henrique.saearp@ime.unicamp.br}{\nolinkurl{henrique.saearp@ime.unicamp.br}} .
\end{center}
\endgroup}
\makeatother
\date{}

\begin{document}

\begin{abstract}
We give a uniform second-order Taylor formula for tensorial Riemannian
$H$-structures in Riemannian normal coordinates whose initial frame is
$H$-adapted.  For a defining multi-tensor $\xi$ with intrinsic torsion $T$,
the linear coefficient is $T_q\diamond\xi$, and the quadratic coefficient
features the covariant derivative of $T$, its iterated diamond action, and
curvature.  The formula applies componentwise to tensors of arbitrary type
and recovers the known $\Gtwo$ and $\Spin(7)$ expansions.  We also give the
almost Hermitian specialisation and prove that vanishing of the complete
second-order torsion contribution forces the intrinsic torsion to vanish
at the centre of the coordinates. 
\end{abstract}

\maketitle
\begingroup\renewcommand{\thefootnote}{\fnsymbol{footnote}}\footnotetext[1]{Corresponding author}\endgroup

\tableofcontents

\newpage

\section{Introduction}

Many Riemannian geometric structures are described by a tensor, or by a finite
collection of tensors, whose pointwise stabiliser is a subgroup
$H\subset\SO(n)$. Familiar examples include positive $3$-forms in dimension
seven, Cayley $4$-forms in dimension eight, and almost Hermitian structures in
even dimensions. Their intrinsic torsion measures the failure of the defining
tensors to be parallel for the Levi--Civita connection.

The purpose of this short note is to express the second-order Taylor expansion
of these tensors in one formula, allowing a defining multi-tensor whose entries
have different tensor types. The linear term records the intrinsic torsion;
the quadratic term involves its covariant derivative, an iterated torsion
action, and the curvature correction arising from the coordinate frame.
The resulting formula also allows us to determine whether the torsion
contributions can be eliminated by imposing a special torsion regime.

We use the intrinsic torsion $T=\nabla^H-\nabla$, where $\nabla^H$ is the
$H$-connection obtained by orthogonally projecting the Levi--Civita connection
onto $\fh$, and the infinitesimal right action $\diamond$ on tensors, so that
$\nabla_X\xi=T_X\diamond\xi$. These conventions are fixed in
Section~\ref{sec:preliminaries}. For a tensor of type $(r,s)$, write
$I=(i_1,\ldots,i_r)$ and $J=(j_1,\ldots,j_s)$ for its upper and lower
multi-indices. Write $u^{\delta^{i_a}_m I}_J$ and $u^I_{\delta^{j_b}_m J}$ for replacement of the indicated
upper or lower entry by $m$, respectively.

\begin{theorem}[Taylor formula]
\label{thm:general-taylor}
Let $(M^n,g,\xi)$ be a smooth Riemannian orientable manifold with tensorial
$H$-structure modelled on a fixed multi-tensor $\xi_0$, with
$H=\Stab_{\SO(n)}(\xi_0)$. Choose normal coordinates
$(x^1,\ldots,x^n)$ centred at $x_0$ whose initial frame is $H$-adapted.
Let $u$ be one entry of the defining multi-tensor $\xi$, of type $(r,s)$.
Its coordinate components have the expansion
\begin{equation}
\label{eq:general-expansion}
 u^I_J(x)
 =u^I_J+(T_q\diamond u)^I_J x^q
   +{}^u\cQ^I_{pqJ}x^px^q+O(|x|^3),
\end{equation}
where the remainder is understood componentwise and every coefficient on the
right-hand side is evaluated at $x_0$. Here
\[
 T_p:=T_{\partial_p},
 \qquad
 \nabla_pT_q:=(\nabla_{\partial_p}T)(\partial_q),
\]
and $\nabla T$ denotes the Levi--Civita covariant derivative after $T$ is viewed
as a section of $T^*M\otimes\fso(TM)$. The quadratic coefficient is
\begin{equation}
\label{eq:general-quadratic}
\begin{split}
 {}^u\cQ^I_{pqJ}
 =&\frac{1}{2}\bigl(\nabla_pT_q\diamond u\bigr)^I_J
   +\frac{1}{2}\bigl(T_q\diamond(T_p\diamond u)\bigr)^I_J\\
 &+\frac16\sum_{b=1}^{s}
   \bigl(R^m_{pqj_b}+R^m_{pj_bq}\bigr)u^I_{\delta^{j_b}_m J}\\
 &-\frac16\sum_{a=1}^{r}
   \bigl(R^{i_a}_{pqm}+R^{i_a}_{pmq}\bigr)u^{\delta^{i_a}_m I}_J.
\end{split}
\end{equation}
Here $R(\partial_p,\partial_q)\partial_b=R^a_{pqb}\partial_a$ and
$R(X,Y)=[\nabla_X,\nabla_Y]-\nabla_{[X,Y]}$. In fact,
$2{}^u\cQ^I_{pqJ}=\partial_p\partial_q u^I_J$ at $x_0$, so this
coefficient is symmetric in $p,q$.
The statement applies componentwise to all entries of $\xi$, irrespective of
their tensor types.
\end{theorem}

The proof is given in
Section~\ref{sec:general-formula}.
It uses only the defining identity for $T$ and the first derivatives of the
Christoffel symbols in normal coordinates. No irreducibility assumption on
the intrinsic-torsion module is needed, and the local statement does not
require $H$ to be connected. In the same section, we examine the vanishing
of the torsion contributions and relate the contracted second-order formula
to harmonic structures and geometric flows.

The $\Gtwo$ expansion was established in \cite[\S 2.6]{DGK25} and was
subsequently used in the construction of the Ricci--harmonic flow in
\cite[\S 3]{D2026}; the $\Spin(7)$ expansion appears in
\cite[\S 7.1]{D2026}. The general tensorial formalism, including the diamond
action and harmonic $H$-structures, is developed in \cite{F2024}.
The extension of the normal-coordinate calculation to other structure groups
is already suggested in \cite[Remark~2.96]{DGK25}. Our contribution is its
componentwise formulation for arbitrary tensorial $H$-structures, together
with the explicit almost Hermitian specialisation below.
Section~\ref{sec:examples} recovers the known $\Gtwo$ and $\Spin(7)$
formulas and derives this specialisation.

\begin{corollary}[$\U(m)$ Taylor formula]
\label{cor:um-taylor-introduction}
Let $(M^{2m},g,J)$ be an almost Hermitian manifold, with fundamental form
$\omega(X,Y)=g(JX,Y)$ and intrinsic torsion
\[
 T_X=-\frac12J\nabla_XJ.
\]
In $\U(m)$-adapted Riemannian normal coordinates centred at $x_0$,
\begin{equation}
\label{eq:metric-expansion}
 g_{ij}(x)=\delta_{ij}+\frac13R_{ipjq}x^px^q+O(|x|^3),
\end{equation}
and
\[
\begin{split}
 \omega_{rs}(x)
 =&\omega_{rs}+(T_q\diamond\omega)_{rs}x^q\\
 &+\frac12\Big[
 (\nabla_pT_q\diamond\omega)_{rs}
 +(T_q\diamond(T_p\diamond\omega))_{rs}\Big]x^px^q\\
 &+\frac16\bigl(
 R^a_{prq}\omega_{as}+R^a_{psq}\omega_{ra}\bigr)x^px^q
 +O(|x|^3),
\end{split}
\]
while
\[
\begin{split}
 J^l_k(x)
 =&J^l_k+(T_q\diamond J)^l_k x^q\\
 &+\frac12\Big[
 (\nabla_pT_q\diamond J)^l_k
 +(T_q\diamond(T_p\diamond J))^l_k\Big]x^px^q\\
 &+\frac16\bigl(
 R^a_{pkq}J^l_a-R^l_{paq}J^a_k\bigr)x^px^q
 +O(|x|^3).
\end{split}
\]
All coefficients on the right-hand sides are evaluated at $x_0$.
\end{corollary}

The derivation is given in the example~\ref{sec:um}, where we also describe the quadratic
torsion action, its relation to the Gray--Hervella classes, and the
Kähler specialisation.
\section{Tensorial \texorpdfstring{$H$}{H}-structures and conventions}
\label{sec:preliminaries}

\subsection{The diamond action}

Let $(V,g_0)$ be an $n$-dimensional Euclidean vector space.  If $u$ is a tensor
of type $(r,s)$ and $a\in\GL(V)$, set
\[
 u\mathbin{\cdot}a
 :=\bigl((a^{-1})^{\otimes r}\otimes(a^*)^{\otimes s}\bigr)u.
\]
For $A\in\End(V)$, the diamond action is the infinitesimal right action
\[
 A\diamond u
 :=\left.\frac{d}{dt}\right|_{t=0}
 u\mathbin{\cdot}\exp(tA).
\]
In components,
\begin{equation}
\label{eq:diamond-components}
\begin{split}
 (A\diamond u)^{i_1\ldots i_r}_{j_1\ldots j_s}
  =&-\sum_{a=1}^{r} A^{i_a}_{m}
  u^{i_1\ldots i_{a-1}m i_{a+1}\ldots i_r}_{j_1\ldots j_s} \\
  &+\sum_{b=1}^{s} A^{m}_{j_b}
  u^{i_1\ldots i_r}_{j_1\ldots j_{b-1}m j_{b+1}\ldots j_s}.
\end{split}
\end{equation}
For a covariant tensor this convention reads
\[
 (A\diamond u)(X_1,\ldots,X_s)
 =\sum_{b=1}^{s}u(X_1,\ldots,AX_b,\ldots,X_s).
\]
For an endomorphism $J$, regarded as a tensor of type $(1,1)$, it gives
\begin{equation}
\label{eq:diamond-J}
  A\diamond J=-[A,J].
\end{equation}
Since this is a right action, its infinitesimal operators satisfy
\begin{equation}
\label{eq:diamond-commutator}
 [A\diamond,B\diamond]u=-[A,B]\diamond u.
\end{equation}
The action is extended componentwise to a multi-tensor.

Whenever a $2$-form is regarded as a skew-adjoint endomorphism, we use the
metric identification
\begin{equation}
\label{eq:two-form-endomorphism}
 g(A_\beta X,Y)=\beta(X,Y),
 \qquad \beta\in\Lambda^2V^*.
\end{equation}

If $A$ is skew-adjoint, the induced infinitesimal action is skew-adjoint on
every tensor representation:
\begin{equation}
\label{eq:skew-adjoint-diamond}
 \langle A\diamond u,v\rangle+\langle u,A\diamond v\rangle=0.
\end{equation}
In particular, $A\diamond g=0$ and
\begin{equation}
\label{eq:negative-square}
 \langle A\diamond(A\diamond u),u\rangle=-|A\diamond u|^2.
\end{equation}

\subsection{Intrinsic torsion}

Fix a finite-dimensional orthogonal direct sum $W$ of tensor representations of
$\SO(n)$, equipped with the induced tensor-product inner products, and a model
multi-tensor $\xi_0\in W$ with
\[
  H=\Stab_{\SO(n)}(\xi_0).
\]
Then $H$ is closed.  Let $\fh$ be its Lie algebra and let
$\fm=\fh^\perp\subset\fso(n)$ with respect to
$\langle A,B\rangle=-\operatorname{tr}(AB)$.  The infinitesimal stabiliser of
$\xi_0$ is $\fh$; consequently,
\begin{equation}
\label{eq:injective-action}
  \fm\longrightarrow W,\qquad A\longmapsto A\diamond\xi_0,
\end{equation}
is injective.

Let $(M^n,g)$ carry an $H$-reduction $Q\subset\Fr(M,g)$ represented by a smooth
multi-tensor $\xi$ modelled on $\xi_0$.  More precisely, the components of
$\xi_x$ in every frame of $Q_x$ are $\xi_0$, and $Q_x$ consists of exactly the
oriented orthonormal frames with this property.  Write
$\fm_Q=Q\times_H\fm$.

In a local $Q$-frame, decompose the Levi--Civita connection form as
$\omega=\omega_{\fh}+\omega_{\fm}$.  The $\fh$-part defines the 
$H$-connection $\nabla^H$.  Its difference from the Levi--Civita connection is
the intrinsic torsion
\[
  T:=\nabla^H-\nabla=-\omega_{\fm}
  \in\Omega^1(M,\fm_Q),
  \qquad
  \nabla^H_XY=\nabla_XY+T_XY.
\]
Its alternation $T_XY-T_YX$ is the affine torsion of $\nabla^H$; the
tensor $T$ itself is the connection difference.  On every tensor bundle, the induced connections satisfy
$\nabla^H_Xu=\nabla_Xu-T_X\diamond u$.  Since $\nabla^H\xi=0$, it follows that
\begin{equation}
\label{eq:intrinsic-torsion}
  \nabla_X\xi=T_X\diamond\xi
\end{equation}
for every vector field $X$, as in \cite[Lemma~1.17]{F2024}.  Conversely,
the injectivity of
\eqref{eq:injective-action} makes $T$ the unique $\fm_Q$-valued one-form with
this property.  

\subsection{Adapted normal coordinates and curvature}

Fix $x_0\in M$ and an $H$-adapted orthonormal frame
$\mathcal E=(e_1,\ldots,e_n)\in Q_{x_0}$.  The coordinates obtained from
\[
  x=(x^1,\ldots,x^n)\longmapsto
  \exp_{x_0}(x^ie_i)
\]
are called $H$-adapted Riemannian normal coordinates at $x_0$.  Thus the
coordinate frame at $x_0$ is adapted to the reduction, but away from $x_0$ it
need not be an $H$-frame.

Repeated indices are summed unless a sum is displayed explicitly.  We use
\[
 R(\partial_i,\partial_j)\partial_k
 =R^a_{ijk}\partial_a,
 \qquad
 R(X,Y)=\nabla_X\nabla_Y-\nabla_Y\nabla_X-\nabla_{[X,Y]},
\]
and set $R_{ijk\ell}=g_{\ell a}R^a_{ijk}$ and
$\Ric_{jk}=g^{i\ell}R_{ijk\ell}$, in agreement with \cite{F2024}.
With this convention, the first derivative
of the Christoffel symbols in Riemannian normal coordinates is
\begin{equation}
\label{eq:normal-Christoffel}
 \left.\partial_p\Gamma^{a}_{qb}\right|_{x_0}
 =\frac13\left(R^{a}_{pqb}+R^{a}_{pbq}\right).
\end{equation}
Indeed, at $x_0$ one has
$R^a_{pqb}=\partial_p\Gamma^a_{qb}-\partial_q\Gamma^a_{pb}$, whereas: 
$\partial_p\Gamma^a_{qb}+\partial_q\Gamma^a_{bp}
+\partial_b\Gamma^a_{pq}=0$.  Together with the first Bianchi identity,
these relations give \eqref{eq:normal-Christoffel}.  

For multi-indices $I=(i_1,\ldots,i_r)$ and $J=(j_1,\ldots,j_s)$, let
$\delta^{i_a}_m I$ and $\delta^{j_b}_m J$ denote replacement of the indicated entry by $m$.
Put
\[
 K_{pq}{}^a_b:=\left.\partial_p\Gamma^a_{qb}\right|_{x_0}
\]
and define the normal-coordinate curvature operator by
$\cC_{pq}u:=K_{pq}\diamond u$.  Thus
\begin{equation}
\label{eq:curvature-operator}
\begin{split}
 (\cC_{pq}u)^I_J
 =&\frac13\sum_{b=1}^{s}
 \left(R^{m}_{pqj_b}+R^{m}_{pj_bq}\right)u^I_{\delta^{j_b}_m J}\\
 &-\frac13\sum_{a=1}^{r}
 \left(R^{i_a}_{pqm}+R^{i_a}_{pmq}\right)u^{\delta^{i_a}_m I}_J.
\end{split}
\end{equation}
All tensors in \eqref{eq:normal-Christoffel} and
\eqref{eq:curvature-operator} are evaluated at $x_0$.

\section{Proof of the main theorem}
\label{sec:general-formula}

Let $\xi$ be the defining tensor of the $H$-structure. Since $$\nabla(\xi_1,\dots,\xi_m)=T\diamond (\xi_1,\dots,\xi_m)=(T\diamond\xi_1,\dots,T\diamond\xi_m)=(\nabla\xi_1,\dots,\nabla\xi_m),$$ we may assume, without loss of generality, that $\xi$ is a single tensor of type $(r,s)$.
This assumption does not assert that each component determines
the intrinsic torsion separately. Indeed, writing
\[
\fh_i=\ker(\cdot\diamond\xi_i),
\qquad
\fh=\bigcap_i\fh_i,
\qquad
\fm=\fh^\perp,
\]
the intersection $\fm\cap\fh_i$ may be nontrivial, so
$T\diamond\xi_i$ need not determine $T$. For instance, a
$U(m)$-structure is defined by $(g,J)$, but $T\diamond g=0$
even when $T\neq0$. This causes no difficulty: our formulas
concern the full $H$-structure and need not decouple its
defining tensors. In particular, the Taylor formula for the
fundamental form $\omega$ will involve $J$, as in the
$U(m)$ example below.

\subsection{Proof of Theorem~\ref{thm:general-taylor}}

We  compute the first and second coordinate
derivatives of $\xi$ in the normal coordinate frame.

\paragraph{\textbf{Claim 1.}}
Let $(x^1,\ldots,x^n)$ be $H$-adapted normal coordinates centred at
$x_0\in M$, let $T$ denote the intrinsic torsion, and let $\nabla$ be the
Levi--Civita connection.  Then, at $x_0$,
\[
 \partial_q \xi^{i_1\ldots i_r}_{j_1\ldots j_s}
 =(T_q\diamond \xi)^{i_1\ldots i_r}_{j_1\ldots j_s}.
\]

\begin{proof}
Writing $\partial_j=\partial/\partial x^j$,
\begin{align*}
    (\nabla_q \xi)^{i_1\ldots i_r}_{j_1\ldots j_s}
    =&\; (\nabla_{\partial_q}\xi)
    (\partial_{j_1}, \ldots,\partial_{j_s}, dx^{i_1}, \ldots,dx^{i_r})\\
    =&\; \partial_q \xi^{i_1\ldots i_r}_{j_1\ldots j_s}\\
    &-\sum_{b=1}^{s}\xi(\partial_{j_1},\ldots,
    \nabla_{\partial_q} \partial_{j_b},\ldots,\partial_{j_s},  dx^{i_1},\ldots,dx^{i_r})\\
    &-\sum_{a=1}^{r}\xi(\partial_{j_1},\ldots,\partial_{j_s},
    dx^{i_1}, \ldots,\nabla_{\partial_q}dx^{i_a},\ldots,dx^{i_r}).
\end{align*}
Since $\nabla_{\partial_q}\partial_{j_b}
=\Gamma^m_{qj_b}\partial_m$ and
$\nabla_{\partial_q}dx^{i_a}=-\Gamma^{i_a}_{qm}dx^m$, this gives
\begin{equation}
\label{eq:first-coordinate-derivative}
 (\nabla_q \xi)^I_J
 =\partial_q \xi^I_J
   -\sum_{b=1}^{s}\Gamma^m_{qj_b}\xi^I_{\delta^{j_b}_m J}
   +\sum_{a=1}^{r}\Gamma^{i_a}_{qm}\xi^{\delta^{i_a}_m I}_J.
\end{equation}
Moreover, at $x_0$ all Christoffel symbols vanish,
because the coordinates are normal. Then, by \eqref{eq:intrinsic-torsion}, 
\[
    \partial_q \xi^I_J
    =(\nabla_q \xi)^I_J
    =(T_q\diamond \xi)^I_J.\qedhere
\]
\end{proof}

\paragraph{\textbf{Claim 2.}}
Under the assumptions of Claim~1, at $x_0$ one has
\begin{equation}
\label{eq:second-coordinate-components}
\begin{split}
    \partial_p\partial_q \xi^{i_1\ldots i_r}_{j_1\ldots j_s}
    =&\bigl(\nabla_pT_q\diamond \xi\bigr)^{i_1\ldots i_r}_{j_1\ldots j_s}
    +\bigl(T_q\diamond(T_p\diamond \xi)\bigr)^{i_1\ldots i_r}_{j_1\ldots j_s}\\
    &+\frac13\sum_{b=1}^{s}
    \bigl(R^m_{pqj_b} +R^m_{pj_bq}\bigr)
    \xi^{i_1\ldots i_r}_{j_1\ldots j_{b-1}m j_{b+1}\ldots j_s}\\
    &-\frac13\sum_{a=1}^{r}
    \bigl(R^{i_a}_{pqm} +R^{i_a}_{pmq}\bigr)
    \xi^{i_1\ldots i_{a-1}m i_{a+1}\ldots i_r}_{j_1\ldots j_s}.
\end{split}
\end{equation}
We denote this expression by $2{}^\xi\cQ^I_{pqJ}$.

\begin{proof}
Differentiating \eqref{eq:first-coordinate-derivative} gives
$$
    \partial_p\partial_q \xi^I_J
    = \partial_p(\nabla_q \xi)^I_J
    +\sum_{b=1}^{s}(\partial_p\Gamma^m_{qj_b})\xi^I_{\delta^{j_b}_m J}
    -\sum_{a=1}^{r}(\partial_p\Gamma^{i_a}_{qm})\xi^{\delta^{i_a}_m I}_J +(\text{terms with } \Gamma).
$$
The terms with $\Gamma$ vanish at $x_0$.  As in Claim~1, at $x_0$,  
$\partial_p(\nabla_q \xi)^I_J
=(\nabla_{\partial_p}(\nabla_{\partial_q}\xi))^I_J$.  Furthermore,
$\nabla_{p}(T_{q})
=(\nabla_pT)_q+\Gamma^m_{pq}T_m$, so compatibility of covariant
differentiation with the diamond action gives, at $x_0$,
\begin{align*}
    \partial_p(\nabla_q \xi)^I_J
    &=\bigl(\nabla_{\partial_p}(T_q\diamond \xi)\bigr)^I_J\\
    &=\bigl((\nabla_pT)_q\diamond \xi\bigr)^I_J
    +\bigl(T_q\diamond\nabla_p\xi\bigr)^I_J\\
    &=\bigl((\nabla_pT)_q\diamond \xi\bigr)^I_J
    +\bigl(T_q\diamond(T_p\diamond \xi)\bigr)^I_J.
\end{align*}
From now on we will write $\nabla_pT_q$ for $(\nabla_p T)_q$. Finally, by \eqref{eq:normal-Christoffel},
\[
 \partial_p\Gamma^m_{qj_b}
 =\frac13\bigl(R^m_{pqj_b}+R^m_{pj_bq}\bigr),
 \qquad
 \partial_p\Gamma^{i_a}_{qm}
 =\frac13\bigl(R^{i_a}_{pqm}+R^{i_a}_{pmq}\bigr)
\]
at $x_0$.  Substitution yields \eqref{eq:second-coordinate-components}.
Equivalently, using \eqref{eq:curvature-operator}, the calculation reads
\begin{equation}
\label{eq:second-coordinate-derivative}
    \left.\partial_p\partial_q\xi^I_J\right|_{x_0}
    =\bigl(\nabla_p\nabla_q \xi+\cC_{pq}\xi\bigr)^I_J,
\end{equation}
with
\begin{equation}
\label{eq:unsymmetrised-hessian}
 \left.\partial_p\partial_q\xi^I_J\right|_{x_0}
 =\bigl(\nabla_pT_q\diamond \xi
   +T_q\diamond(T_p\diamond \xi)+\cC_{pq}\xi\bigr)^I_J.
 \qedhere
\end{equation}
\end{proof}

The ordinary multivariable Taylor formula is
\[
 \xi^I_J(x)=\xi^I_J+\partial_q\xi^I_Jx^q
 +\frac12\partial_p\partial_q\xi^I_Jx^px^q+O(|x|^3).
\]
Substituting Claims~1 and~2 proves
\eqref{eq:general-expansion}--\eqref{eq:general-quadratic}, and hence
Theorem~\ref{thm:general-taylor}.  Equality of mixed partial derivatives 
also shows directly that ${}^\xi\cQ^I_{pqJ}={}^\xi\cQ^I_{qpJ}$.

\begin{remark}[Expanded curvature term]
\label{rem:expanded-curvature}
Since $x^px^q$ is symmetric and $R^a_{pq b}$ is skew in $p,q$, the expansion
may equivalently be written as
\begin{align}
 \xi^I_J(x)=&\xi^I_J+(T_q\diamond \xi)^I_J x^q \notag\\
 &+\frac12\bigl(\nabla_pT_q\diamond \xi
       +T_q\diamond(T_p\diamond \xi)\bigr)^I_J x^px^q \notag\\
 &+\frac16\left[
   \sum_{b=1}^{s}R^m_{p j_b q}\xi^I_{\delta^{j_b}_m J}
   -\sum_{a=1}^{r}R^{i_a}_{p m q}\xi^{\delta^{i_a}_m I}_J
  \right]x^px^q+O(|x|^3).
\label{eq:expanded-general-formula}
\end{align}
\end{remark}

\subsection{Vanishing of the torsion terms}

We now ask whether the Taylor formula simplifies formally when some of the torsion terms in \eqref{eq:expanded-general-formula} vanish. Following \cite{F2024},
write $\Delta=-\nabla^*\nabla$ for the rough Laplacian.
Contracting \eqref{eq:unsymmetrised-hessian} at $x_0$ gives
\[
 \sum_p\partial_p\partial_p\xi
 =\Delta\xi-\frac13\Ric^\sharp\diamond\xi,
 \qquad
 \Delta\xi
 =(\operatorname{div}T)\diamond\xi
  +\sum_pT_p\diamond(T_p\diamond\xi).
\]
Thus, vanishing of the second-order torsion contribution at $x_0$ 
implies $\Delta\xi(x_0)=0$. The endomorphism $\operatorname{div}T$ is
skew-adjoint, so \eqref{eq:skew-adjoint-diamond} gives
$\langle(\operatorname{div}T)\diamond\xi,\xi\rangle=0$.
Consequently, \eqref{eq:negative-square} yields 
\[
 \langle\Delta\xi,\xi\rangle
 =\sum_p\langle T_p\diamond(T_p\diamond\xi),\xi\rangle
 =-\sum_p|T_p\diamond\xi|^2.
\]
The injectivity of \eqref{eq:injective-action} then forces
$T(x_0)=0$. In particular, $\operatorname{div}T=0$ or
$\nabla T=0$ alone do not remove the quadratic torsion
contribution when $T\neq0$.

Conversely, $T(x_0)=0$ removes the linear term and the iterated
quadratic $T$-term, but the term involving
$\nabla T(x_0)$ may remain. If $T$ vanishes on a
neighbourhood of $x_0$, then only the curvature terms in the Taylor formula remain.

\subsection{Harmonic structures and geometric flows}
\label{sec:applications}

The linear Taylor coefficient of the complete defining multi-tensor vanishes
in every direction precisely when $T(x_0)=0$.  The trace of the second order term is the rough-Laplacian identity of
\cite[Lemma~1.21]{F2024} together with the associated Ricci curvature correction.  Although $\nabla$ need not
preserve $\fm_Q$, its contraction is $\fm_Q$-valued: the connection-difference
identity
$(\nabla_XT)_Y=(\nabla^H_XT)_Y-[T_X,T_Y]+T_{T_XY}$ implies
\begin{equation}
\label{eq:divergence-in-m}
 \operatorname{div}T
 =\sum_p(\nabla^H_{e_p}T)_{e_p}
   +T_{\sum_pT_{e_p}e_p}\in\fm_Q.
\end{equation}

For the variational interpretation, assume that $M$ is closed, $n>2$, and
$H$ is connected, as in \cite{F2024}.  For isometric variations with
$\left.\partial_t\xi(t)\right|_{t=0}=C\diamond\xi$, where
$C\in\Gamma(\fm_Q)$, the intrinsic-torsion energy satisfies
\[
 \mathcal E(\xi)=\frac12\int_M|T|^2\vol_g,
 \qquad
 \left.\frac{d}{dt}\right|_{t=0}\mathcal E(\xi(t))
 =-\int_M\langle\operatorname{div}T,C\rangle\vol_g;
\]
see \cite[Corollary~1.46]{F2024}.  Thus harmonicity is characterised by $\operatorname{div}T=0$,
which alone does not eliminate the quadratic torsion contribution.
Combining $\operatorname{div}T$ with the Ricci action gives the
Ricci--harmonic evolution discussed in \cite{FFLMSE26}. In our general
setting the metric is retained explicitly, so this is an evolution of
the compatible pair $(g,\xi)$:
\[
 \partial_tg=-2\Ric,
 \qquad
 \partial_t\xi=(-\Ric^\sharp+\operatorname{div}T)\diamond\xi.
\]
Equivalently, one may include $g$ in the defining multi-tensor; the metric
equation then follows from $A\diamond g=2g(A_{\mathrm{sym}}\cdot,\cdot)$.
The condition $H=\Stab_{\SO(n)}(\xi_0)$ alone does not assert that $\xi$
determines the metric.  The Taylor expansions motivate related flows in
\cite{D2026}, but do not by themselves specify a flow on the space of
$H$-structures.  In particular, an iterated diamond action need not be a
single diamond action, as the $\Spin(7)$ and $\U(m)$ examples below show. Thus an endomorphism-valued quadratic expression $B(T)$ satisfying
\[
 B(T)\diamond\xi=\sum_pT_p\diamond(T_p\diamond\xi)
\]
need not exist. The complete traced iterated action cannot, in general,
be inserted into such a flow as a single diamond action.
Natural lower-order velocities of the form $B(\xi,T)\diamond\xi$, quadratic
in $T$, are treated in \cite{FFLMSE26}.

\section{Examples}
\label{sec:examples}

\paragraph{\textbf{\texorpdfstring{$\Gtwo$}{G2}-structures.}}
\label{sec:g2}

A $\Gtwo$-structure on a $7$-manifold is represented by a positive $3$-form
$\varphi$.  It determines a metric and orientation; with the convention used
here,
\[
 (X\mathbin{\lrcorner}\varphi)\wedge
 (Y\mathbin{\lrcorner}\varphi)\wedge\varphi
 =-6g(X,Y)\vol_g.
\]
Writing $\psi=*\varphi$, the $\Gtwo$-decomposition of $2$-forms is
\[
 \Lambda^2=\Lambda^2_7\oplus\Lambda^2_{14},
 \qquad
 \Lambda^2_7=\{\beta:*(\beta\wedge\varphi)=-2\beta\},
 \quad
 \Lambda^2_{14}=\{\beta:*(\beta\wedge\varphi)=\beta\}.
\]
Here $\Lambda^2_{14}\simeq\mathfrak{g}_2$ and
$\Lambda^2_7\simeq\mathfrak{g}_2^\perp$.  Thus the intrinsic torsion is
\[
 T\in\Omega^1(M,\Lambda^2_7).
\]

\begin{corollary}[$\Gtwo$ Taylor formula]
\label{cor:g2-taylor}
In $\Gtwo$-adapted Riemannian normal coordinates at $x_0$,
\begin{equation}
\label{eq:g2-expansion}
\begin{split}
 \varphi_{ijk}(x)
 =&\varphi_{ijk}+(T_q\diamond\varphi)_{ijk}x^q\\
 &+\frac12\Big[
 (\nabla_pT_q\diamond\varphi)_{ijk}
 +(T_q\diamond(T_p\diamond\varphi))_{ijk}\Big]x^px^q\\
 &+\frac16\bigl(
 R^a_{piq}\varphi_{ajk}
 +R^a_{pjq}\varphi_{iak}
 +R^a_{pkq}\varphi_{ija}\bigr)x^px^q
 +O(|x|^3).
\end{split}
\end{equation}
All coefficients on the right-hand side are evaluated at $x_0$.
\end{corollary}

The normal-coordinate expansion of Dwivedi, Gianniotis and
Karigiannis \cite[Theorem~2.91]{DGK25} is recovered exactly
from our general formula after translating the torsion
conventions. To make this comparison explicit, denote their
full torsion tensor by $\mathcal T$, so that
\[
 \nabla_q\varphi_{ijk}=\mathcal T_{qm}\psi_{mijk}.
\]
With $T_{q;ij}:=g(T_qe_i,e_j)$, the conversion is
\begin{equation}
\label{eq:g2-full-torsion}
 T_{q;ij}=-\frac13\mathcal T_{qm}\varphi_{mij},
 \qquad
 (T_q\diamond\varphi)_{ijk}=\mathcal T_{qm}\psi_{mijk},
 \qquad
 \mathcal T_{qm}=\frac1{24}(\nabla_q\varphi_{ijk})\psi_{mijk};
\end{equation}
see also \cite[Example~1.19]{F2024}.
Since $\nabla$ commutes with the Hodge star,
\[
 \nabla_p\psi_{mijk}
 =-\mathcal T_{pm}\varphi_{ijk}
  +\mathcal T_{pi}\varphi_{mjk}
  -\mathcal T_{pj}\varphi_{mik}
  +\mathcal T_{pk}\varphi_{mij}.
\]
Substituting these identities into the general Taylor formula
gives
\begin{equation}
\label{eq:g2-full-expansion}
 \varphi_{ijk}(x)
 =\varphi_{ijk}+\mathcal T_{qm}\psi_{mijk}x^q
  +\mathcal P_{pqijk}x^px^q+O(|x|^3),
\end{equation}
where
\begin{align}
 \mathcal P_{pqijk}={}&
 \frac12(\nabla_p\mathcal T_{qm})\psi_{mijk}
 -\frac12(\mathcal T\mathcal T^{\top})_{pq}\varphi_{ijk}
 \notag\\
 &+\frac12\mathcal T_{pm}\bigl(
   \mathcal T_{qi}\varphi_{mjk}
  +\mathcal T_{qj}\varphi_{mki}
  +\mathcal T_{qk}\varphi_{mij}\bigr)
 \notag\\
 &+\frac16\bigl(
 R_{piqm}\varphi_{mjk}
 +R_{pjqm}\varphi_{mki}
 +R_{pkqm}\varphi_{mij}\bigr),
\label{eq:g2-full-quadratic}
\end{align}
and
\[
 (\mathcal T\mathcal T^{\top})_{pq}
 =\mathcal T_{pm}\mathcal T_{qm}.
\]
All coefficient tensors are evaluated at the centre of the
normal coordinates. Only the symmetric part in $p,q$
contributes to the expansion. Thus
\eqref{eq:g2-full-expansion}--\eqref{eq:g2-full-quadratic}
coincide with \cite[(2.92)--(2.93)]{DGK25}.
Although $\Lambda^2_7$ is irreducible, the full intrinsic-torsion module is not:
\begin{equation}
\label{eq:g2-torsion-decomposition}
 T^*M\otimes\Lambda^2_7
 \simeq \mathcal W_1\oplus\mathcal W_2\oplus\mathcal W_3\oplus\mathcal W_4
 \simeq \Lambda^0\oplus\Lambda^2_{14}\oplus\Lambda^3_{27}\oplus\Lambda^1.
\end{equation}
Equivalently, the four torsion forms
$\tau_0\in\Omega^0$, $\tau_1\in\Omega^1$,
$\tau_2\in\Omega^2_{14}$, and $\tau_3\in\Omega^3_{27}$ are characterised, in a
standard normalisation, by
\begin{equation}
\label{eq:g2-torsion-forms}
 d\varphi=\tau_0\psi+3\tau_1\wedge\varphi+*\tau_3,
 \qquad
 d\psi=4\tau_1\wedge\psi+\tau_2\wedge\varphi.
\end{equation}
For further reading on the subject we refer to the classical literature by Fernández-Gray and Karigiannis \cite{FG82,Kar09}, respectively.  These forms determine $T$, and hence may be substituted
in \eqref{eq:g2-expansion}.  The symmetric iterated term lies in
$\Lambda^3_1\oplus\Lambda^3_{27}$: this follows from equivariance and
$\operatorname{Sym}^2(\Lambda^2_7)\simeq\Lambda^0\oplus\Lambda^3_{27}$.
This algebraic restriction does not imply its vanishing when $T\ne0$, by
\eqref{eq:negative-square}.
In particular, $T=0$ if and only if all four torsion forms vanish, equivalently
if and only if $d\varphi=d\psi=0$.

\paragraph{\textbf{\texorpdfstring{$\Spin(7)$}{Spin(7)}-structures.}}
\label{sec:spin7}

A $\Spin(7)$-structure on an $8$-manifold is represented by a Cayley $4$-form
$\Phi$.  The form determines a metric and orientation and is self-dual:
$*\Phi=\Phi$.  The relevant decomposition is
\[
 \Lambda^2=\Lambda^2_7\oplus\Lambda^2_{21},
 \qquad
\Lambda^2_7=\{\beta:*(\beta\wedge\Phi)=-3\beta\},
 \quad
\Lambda^2_{21}=\{\beta:*(\beta\wedge\Phi)=\beta\}.
\]
These eigenvalues fix the overall sign convention for the invariant Cayley
tensor.
Here $\Lambda^2_{21}\simeq\mathfrak{spin}(7)$ and
$\Lambda^2_7\simeq\mathfrak{spin}(7)^\perp$, so
$T\in\Omega^1(M,\Lambda^2_7)$.
Writing $T_{q;ij}:=g(T_qe_i,e_j)$, the defining relation and its inverse are
\begin{align}
 \nabla_q\Phi_{ijkl}
  =&(T_q\diamond\Phi)_{ijkl} \notag\\
  =&T_{q;ia}\Phi_{ajkl}+T_{q;ja}\Phi_{iakl}
      +T_{q;ka}\Phi_{ijal}+T_{q;la}\Phi_{ijka},
\label{eq:spin7-torsion}\\
 T_{q;ab}=&\frac1{96}(\nabla_q\Phi_{ajkl})\Phi_{bjkl}.
\label{eq:spin7-torsion-inverse}
\end{align}
Indeed, for every $A\in\Lambda^2_7$,
\[
 \bigl((A\diamond\Phi)_{ajkl}\bigr)\Phi_{bjkl}=96A_{ab};
\]
see also \cite[(2.14)--(2.15)]{D2026}.

\begin{corollary}[$\Spin(7)$ Taylor formula]
\label{cor:spin7-taylor}
In $\Spin(7)$-adapted Riemannian normal coordinates at $x_0$,
\begin{equation}
\label{eq:spin7-expansion}
\begin{split}
 \Phi_{ijkl}(x)
 =&\Phi_{ijkl}+(T_q\diamond\Phi)_{ijkl}x^q\\
 &+\frac12\Big[
 (\nabla_pT_q\diamond\Phi)_{ijkl}
 +(T_q\diamond(T_p\diamond\Phi))_{ijkl}\Big]x^px^q\\
 &+\frac16\bigl(
 R^a_{piq}\Phi_{ajkl}
 +R^a_{pjq}\Phi_{iakl}
 +R^a_{pkq}\Phi_{ijal}
 +R^a_{plq}\Phi_{ijka}\bigr)x^px^q
 +O(|x|^3).
\end{split}
\end{equation}
All coefficients on the right-hand side are evaluated at $x_0$.
\end{corollary}

This coincides with \cite[Theorem~7.1]{D2026}.

The intrinsic-torsion module has two irreducible summands:
\begin{equation}
\label{eq:spin7-torsion-decomposition}
 T^*M\otimes\Lambda^2_7\simeq\Lambda^1\oplus\Lambda^3_{48},
 \qquad
 \Lambda^3_{48}=\{\gamma\in\Lambda^3:\gamma\wedge\Phi=0\}.
\end{equation}
Equivalently, one may define a Lee form $\vartheta$ and a form
$\tau_3\in\Omega^3_{48}$ by decomposing
\begin{equation}
\label{eq:spin7-torsion-forms}
 d\Phi=\vartheta\wedge\Phi+*\tau_3,
\end{equation}
with the normalisation of $\vartheta$ fixed by this identity.  See
\cite{F86}.  A $\Spin(7)$-structure is torsion-free if and only if $d\Phi=0$,
equivalently if and only if $\vartheta=\tau_3=0$.

For $A,B\in\Lambda^2_7$, equivariance and the decompositions
\[
 \Lambda^2_7\otimes\Lambda^2_7\simeq\Lambda^0\oplus\Lambda^2_{21}\oplus\Lambda^4_{27}, \qquad
\Lambda^4=\Lambda^4_1\oplus\Lambda^4_7\oplus\Lambda^4_{27}\oplus\Lambda^4_{35}
\]
show that $A\diamond(B\diamond\Phi)$ lies in
$\Lambda^4_1\oplus\Lambda^4_{27}$.  In contrast, the image of
$\End(TM)\ni A\mapsto A\diamond\Phi$ is
$\Lambda^4_1\oplus\Lambda^4_7\oplus\Lambda^4_{35}$; see
\cite[Proposition~2.2]{D2026}.  Thus the quadratic torsion term cannot in
general be replaced by a single endomorphism acting on $\Phi$.
To see that the $\Lambda^4_{27}$ component really occurs, use the
standard $\SU(4)$ model on $\mathbb R^8$ with K\"ahler form $\omega_0$,
complex structure $J_0$ and complex volume form $\Omega_0$, writing
\[
 \Phi=-\frac12\omega_0^2-\operatorname{Re}\Omega_0
\]
in the sign convention above. Since
$*(\omega_0\wedge\Phi)=-3\omega_0$, the endomorphism $J_0$
corresponds to an element of $\Lambda^2_7$. Its right action fixes
$\omega_0$ and sends $\Omega_0$ to $e^{4it}\Omega_0$, whence
\[
 J_0\diamond(J_0\diamond\Phi)=16\operatorname{Re}\Omega_0.
\]
This is not a scalar multiple of $\Phi$. As it belongs to
$\Lambda^4_1\oplus\Lambda^4_{27}$, its $\Lambda^4_{27}$ component
is nonzero, proving the asserted obstruction.

\paragraph{\textbf{\texorpdfstring{$\U(m)$}{U(m)}-structures.}}
\label{sec:um}

The case $m=1$ is vacuous from the present viewpoint, since
$\U(1)=\SO(2)$ and $\mathfrak u(1)^\perp=0$.  We therefore assume $m\geq2$.
A $\U(m)$-structure on a $2m$-manifold is an almost Hermitian structure $(g,J)$,
with $J^2=-\id$, $g(JX,JY)=g(X,Y)$, and fundamental form
\[
 \omega(X,Y)=g(JX,Y).
\]
We retain the redundant defining multi-tensor $(g,J,\omega)$ in order to display
the expansions of all three associated tensors; its entries are always treated
separately because they have different tensor types.  Under
$\Lambda^2\simeq\fso(2m)$,
\[
\begin{aligned}
 \mathfrak{u}(m)&=\{A\in\fso(2m):AJ=JA\},\\
 \fm=\mathfrak{u}(m)^\perp
 &=\{A\in\fso(2m):AJ=-JA\}.
\end{aligned}
\]
The intrinsic torsion takes values in $\fm$.  By
\eqref{eq:diamond-J},
\begin{equation}
\label{eq:um-torsion-J}
 \nabla_XJ=T_X\diamond J=-[T_X,J]=2JT_X,
 \qquad
 T_X=-\frac12J\nabla_XJ=\frac12(\nabla_XJ)J.
\end{equation}
Since $\nabla g=0$, one also has
\begin{equation}
\label{eq:um-torsion-omega}
 (\nabla_X\omega)(Y,Z)=g((\nabla_XJ)Y,Z)
  =(T_X\diamond\omega)(Y,Z).
\end{equation}
For a fixed almost Hermitian structure, these identities recover $T$
from $\nabla J$, or equivalently from $\nabla\omega$ after using
$g$ to identify forms and endomorphisms \cite[Example~1.18]{F2024}.
For a general multi-tensor structure, the derivative of one entry need
not determine the full intrinsic torsion. For example, retaining only
one almost complex structure from the defining data of an $\SU(2)$-structure
enlarges the stabiliser and can lose torsion components.

\begin{corollary}[$\U(m)$ Taylor formula]
\label{cor:um-taylor}
In $\U(m)$-adapted Riemannian normal coordinates at $x_0$, the three defining
tensors have the following componentwise expansions.  The metric expansion is
\eqref{eq:metric-expansion}; the remaining entries are
\begin{equation}
\label{eq:omega-expansion}
\begin{split}
 \omega_{rs}(x)
 =&\omega_{rs}+(T_q\diamond\omega)_{rs}x^q\\
 &+\frac12\Big[
 (\nabla_pT_q\diamond\omega)_{rs}
 +(T_q\diamond(T_p\diamond\omega))_{rs}\Big]x^px^q\\
 &+\frac16\bigl(
 R^a_{prq}\omega_{as}+R^a_{psq}\omega_{ra}\bigr)x^px^q
 +O(|x|^3),
\end{split}
\end{equation}
and
\begin{equation}
\label{eq:J-expansion}
\begin{split}
 J^l_k(x)
 =&J^l_k+(T_q\diamond J)^l_k x^q\\
 &+\frac12\Big[
 (\nabla_pT_q\diamond J)^l_k
 +(T_q\diamond(T_p\diamond J))^l_k\Big]x^px^q\\
 &+\frac16\bigl(
 R^a_{pkq}J^l_a-R^l_{paq}J^a_k\bigr)x^px^q
 +O(|x|^3).
\end{split}
\end{equation}
All coefficients on the right-hand sides are evaluated at $x_0$.
\end{corollary}

The linear terms are simply $\nabla_q\omega$ and $\nabla_qJ$, by
\eqref{eq:um-torsion-J}--\eqref{eq:um-torsion-omega}.  The metric has no torsion
terms because it is parallel for the Levi--Civita connection.

For $A,B\in\mathfrak u(m)^\perp$, the relations $AJ=-JA$ and $BJ=-JB$
give
\begin{equation}
\label{eq:um-double-action}
 A\diamond(B\diamond J)=[A,[B,J]]=2J(AB+BA)\in\mathfrak u(m).
\end{equation}
Under the metric identification, the same statement holds for the double
action on $\omega$.  In particular,
$\sum_pT_p\diamond(T_p\diamond J)=4J\sum_pT_p^2$ commutes with $J$,
whereas every single action $A\diamond J$ anticommutes with $J$.
Consequently this sum is a single diamond action on $J$ only if it vanishes,
which forces $T=0$; compare \cite[Example~1.22]{F2024}.

The intrinsic-torsion module of an almost Hermitian structure has the
Gray--Hervella decomposition
\begin{equation}
\label{eq:gray-hervella}
 T^*M\otimes\mathfrak{u}(m)^\perp
 =\mathcal W_1\oplus\mathcal W_2\oplus\mathcal W_3\oplus\mathcal W_4.
\end{equation}
For $m\geq3$, the summands are irreducible real $\U(m)$-modules.  In the usual
notation,
\[
\begin{aligned}
 \mathcal W_1&\simeq
 \{\gamma\in\Lambda^3T^*M:\gamma_{\mathbb C}
   \in\Lambda^{3,0}\oplus\Lambda^{0,3}\},\\
 \mathcal W_3&\simeq
 \{\gamma\in\Lambda^3T^*M:\gamma_{\mathbb C}
   \in\Lambda^{2,1}_0\oplus\Lambda^{1,2}_0\},\\
 \mathcal W_4&\simeq T^*M.
\end{aligned}
\]
Here $\gamma_{\mathbb C}$ denotes complexification and the subscript $0$
denotes the primitive component.
The remaining component has a convenient intrinsic description.  For
$S\in T^*M\otimes\mathfrak{u}(m)^\perp$, set
\[
 s(X,Y,Z)=g(S_XY,Z),
 \qquad
 \operatorname{Alt}(S)(X,Y,Z)
 =\frac13\sum_{\mathrm{cyc}}s(X,Y,Z).
\]
Then $\mathcal W_2=\ker(\operatorname{Alt})$. For a detailed description of these classes, we refer to Gray--Hervella's paper \cite{GH80}.  The real
dimensions are
\[
\begin{aligned}
 \dim_{\mathbb R}\mathcal W_1&=\frac{m(m-1)(m-2)}3,
 &\qquad
 \dim_{\mathbb R}\mathcal W_2&=\frac{2m(m^2-1)}3,\\
 \dim_{\mathbb R}\mathcal W_3&=m(m+1)(m-2),
 &
 \dim_{\mathbb R}\mathcal W_4&=2m.
\end{aligned}
\]
For $m\geq3$, structures of pure classes
$\mathcal W_1,\mathcal W_2,\mathcal W_3,\mathcal W_4$ are, respectively,
nearly K\"ahler, almost K\"ahler, balanced Hermitian, and locally conformally
K\"ahler.  We normalise the Lee form $\theta$ by
\[
 d\omega^{m-1}=\theta\wedge\omega^{m-1}.
\]
In complex dimension two, $\mathcal W_1=\mathcal W_3=0$, and $J$ is integrable
precisely when the $\mathcal W_2$-component vanishes.  Thus every Hermitian
surface is of class $\mathcal W_4$; it is locally conformally K\"ahler precisely
when $d\theta=0$.  In every dimension the structure is K\"ahler precisely when
$T=0$.

If $(g,J,\omega)$ is K\"ahler, then $T=0$, and
\eqref{eq:omega-expansion} reduces to
\begin{equation}
\label{eq:kahler-osculation}
\begin{split}
 \omega_{rs}(x)=\omega_{rs}+\frac16\big[&
 R^a_{prq}\omega_{as}
 +R^a_{psq}\omega_{ra}\big]x^px^q
 +O(|x|^3).
\end{split}
\end{equation}
Thus the K\"ahler form agrees with its flat model to first order; its first
possible non-constant term is quadratic.  The coordinates in \eqref{eq:kahler-osculation} are real Riemannian normal
coordinates. Holomorphic normal coordinates give a related
curvature expansion; see \cite{HU05}. These are different coordinate
normalisations, so their quadratic component formulas should not be
identified without the corresponding coordinate transformation.

\section*{Acknowledgements}
\par\medskip\noindent\textbf{Funding.}\quad
LD is supported by the São Paulo Research Foundation (FAPESP)
[grant number 2025/10187-8]. HSE is supported by FAPESP
[grant numbers 2020/09838-0, 2021/04065-6, 2024/00923-6] and by the
Brazilian National Council for Scientific and Technological Development
(CNPq) [grant number 307145/2025-5]. HSE is Vice-Director of the Centro Brasileiro de Geometria (CBG).
The funders had no role in the research, preparation of the manuscript or decision to submit the article for publication.
\par\medskip\noindent\textbf{Data availability.}\quad
No datasets were generated or analysed in this study. The results
are established by mathematical arguments presented in the article.

\par\medskip\noindent\textbf{Declaration of competing interests.}\quad
The authors declare that they have no known competing financial interests or personal relationships that could have appeared to influence the work reported in this paper.

\par\medskip\noindent\textbf{Declaration of generative AI and AI-assisted technologies in the manuscript preparation process.}\quad
ChatGPT (OpenAI) was used to assist with \LaTeX{} organisation and
formatting, language editing, notation and convention checks, and a
mathematical consistency review. This assistance included suggestions
for clarifying the torsion-vanishing argument, distinguishing coordinate
normalisations and metric data, and developing the explicit algebraic
example in the $\Spin(7)$ discussion. Selected algebraic identities were
also checked by exact computer calculations. AI-generated suggestions
and computational checks do not constitute independent mathematical
verification. Responsibility for checking the arguments, references
and final wording, and for the decision to disseminate the work, rests
entirely with the authors.

\par\medskip\noindent\textbf{Other acknowledgements.}\quad
We thank Udhav Fowdar for useful discussions.

\bibliographystyle{unsrtnat}
\bibliography{references}

\end{document}